\magnification=\magstep1
\def\mapname#1{\ \smash{\mathop{\longrightarrow}\limits^{#1}}\ } 
\overfullrule=0pt

\def\res{{\rm Res}} 
\def\vol{{\rm vol}} 
\def\P{{\bf P}}     
\def\C{{\bf C}}     
\def\R{{\bf R}}     
\def\Q{{\bf Q}}     
\def\Z{{\bf Z}}     
\def\N{{\bf N}}     
\def\n{{\eta}}     
\def\ox{{\cal O}_X} 
\def\LL{{\cal L}} 
\def\RR{{\cal R}}
\def\AA{{\cal A}} 
\def\noi{\noindent}

\font\bigbf=cmbx12 scaled \magstep1  
\def\linebreak{\relax\ifhmode\unskip\break\else
    \nonhmodeerr@\linebreak\fi}
\def\newline{\relax\ifhmode\null\hfil\break
\else\nonhmodeerr@\newline\fi}

\def\frac#1#2{{#1\over#2}}
\def\ip#1#2{\langle #1, #2\rangle}
\def\Efes{F_0,\ldots,F_n}

\def\bi{b_i}

\def\res{{\rm Res}}

\centerline {\bigbf Residues and Resultants} 

\vskip 1cm

\centerline{{\bf  Eduardo Cattani\rm \footnote{${}^1 \!\!$}
{Supported in part by 
the N.S.F.}}}
\centerline{ Department of Mathematics and Statistics}
\centerline{ University of Massachusetts}
\centerline{ Amherst, MA 01003, \ U.S.A. }
\centerline{\tt cattani@math.umass.edu }

\vskip .5cm

\centerline{\bf Alicia Dickenstein\rm \footnote{${}^2 \!\!$}
{Supported in part by 
UBACYT and CONICET, Argentina.}}
\centerline{Departamento de Matem\'atica, F.C.E. y N. }
\centerline{Universidad de Buenos Aires}
\centerline{ Ciudad Universitaria - Pabell\'on I } 
\centerline{(1428) Buenos Aires, \ Argentina }
\centerline{\tt alidick@dm.uba.ar}

\vskip .5cm

\centerline{\bf Bernd Sturmfels\rm \footnote{${}^3 \!\!$}
{Supported in part by 
the N.S.F.~and the David and Lucile Packard Foundation.}}
\centerline{ Department of Mathematics}
\centerline{ University of California }
\centerline{ Berkeley, CA 94720, \  U.S.A. }
\centerline{\tt bernd@math.berkeley.edu }

\vskip 2cm

\centerline{\bf Abstract}
\medskip
\noi { \sl Resultants, Jacobians and residues are basic invariants of 
multivariate polynomial systems. We examine their interrelations in
the context of toric geometry. The global residue in the torus, 
studied
by Khovanskii, is the sum over local Grothendieck residues at the
zeros of $n$ Laurent polynomials in $n$ variables. Cox 
introduced the
related notion of the toric residue relative to
$n+1$ divisors on an $n$-dimensional toric variety.
We establish denominator formulas in terms of sparse resultants for
both the toric residue and the global residue in the torus. A
byproduct is a determinantal formula for resultants based on 
Jacobians. \/}

\vfill
\eject

\noi {\bf \S 0. Introduction }
\smallskip
\noindent
Resultants, Jacobians and residues are fundamental
invariants associated  with systems of multivariate polynomial
equations. We shall investigate relationships among
these three invariants in the context of toric geometry.
The study of {\it global residues in the torus} has its origin 
in the work
of Khovanskii [K2]. The global residue is the sum over local
Grothendieck residues at the common roots of
 $n$ Laurent polynomials in $n$ variables; see (3.8) and (3.10).
The related notion of the {\it toric residue} was introduced
by Cox [C2] and subsequently studied in [CCD]. The toric residue
is associated with $n+1$  divisors on an $n$-dimensional
projective toric variety. For our purposes here
it suffices to consider divisors that
are multiples of a fixed ample divisor $\beta$.
An algorithmic link between these two notions of residue
(``toric'' versus ``in the torus'') was established in [CD].

The main results of this paper are denominator formulas
for toric residues (Theorem~1.4) and for 
residues in the torus (Theorem~3.2). In each case the
denominator is given in terms 
of  {\it sparse resultants}. These
resultants are naturally associated
with sparse systems of Laurent polynomials,
or with line bundles on toric varieties. They were
introduced by Gel'fand, Kapranov and Zelevinsky [GKZ]
and further studied in [KSZ],[PSt],[S1],[S2].
In \S 4 we present new determinantal formulas
for sparse resultants based on Jacobians.

One general objective of our work is to develop
computational techniques, which may
ultimately enter into the design of algorithms
for solving polynomial equations. Classical results
on residues, Jacobians and resultants are limited to
dense equations, in which case the underlying toric
variety is complex projective $n$-space $\,{\bf P}^n$. In that 
classical
case our denominator formula appeared already in the work of
Ang\'eniol [A] and Jouanolou [J1],[J3].
Our results also extend the work of Gel'fond-Khovanskii [GK]
and Zhang [Z], who studied residues in the torus for the special
case when all facet resultants are monomials.

We illustrate our results for two generic
quadratic equations in two complex variables:
$$ \eqalign{
f_1 \quad = \quad &
 a_0  x^2 \,+\, a_1  x y \,+\, a_2  y^2 \,+\, a_3  x  
\,+\, a_4  y \,+\, a_5\,,
 \cr
f_2 \quad = \quad &
b_0 x^2 \,+\, b_1  x y \,+\, b_2  y^2 \,+\, b_3  x  
\,+\, b_4  y \,+\, b_5\,.
\cr} \eqno (0.1) $$
They have four common zeros $(x_i,y_i), i=1,\ldots,4$, in the
algebraic torus $({\bf C}^*)^2$, and the (affine toric) Jacobian
$\,J^T (x,y) \,:= \,xy(\,
{\partial f_1 \over \partial x} { \partial f_2 \over \partial y} -
{\partial f_1 \over \partial y} { \partial f_2 \over \partial x})\,$
is non-zero at these four points. Consider
any Laurent monomial $x^i y^j$. The {\it global residue} 
is the expression
$$ \res^T_f(x^i y^j) \quad := \quad
{x_1^i y_1^j \over  J^T(x_1,y_1)} \,+\,
{x_2^i y_2^j \over  J^T(x_2,y_2)} \,+\,
{x_3^i y_3^j \over  J^T(x_3,y_3)} \,+\,
{x_4^i y_4^j \over  J^T(x_4,y_4)} \,. \eqno (0.2) $$
This is a rational function in the twelve
indeterminates $a_0,\ldots,a_5,b_0,\ldots,b_5 $.
Theorem~3.2 implies that there
exists a polynomial $P_{ij}(a_0,\ldots,a_5,b_0,\ldots,b_5 ) $ 
such that
(0.2) equals
$$ { P_{ij} \over  \RR_\infty^{\max\{0,i+j-3\}} \cdot 
\RR_x^{\max\{0,1-i\}}
 \cdot \RR_y^{\max\{0,1-j\} }} \, , \eqno (0.2') $$
where the prime divisors in the denominator are the
{\it facet resultants}
$$  \eqalign{ \RR_\infty \quad & = \quad
  a_0^2 b_2^2 - a_0 a_1 b_1 b_2
- 2 a_0 a_2 b_0 b_2
+ a_0 a_2 b_1^2
+ a_1^2 b_0 b_2
- a_1 a_2 b_0 b_1
+ a_2^2 b_0^2 \, , \cr
\RR_x  \quad & = \quad
a_0^2 b_5^2 - a_0 a_3 b_3 b_5
 - 2 a_0 a_5 b_0 b_5 + a_0 a_5 b_3^2
+ a_3^2 b_0 b_5 - a_3 a_5 b_0 b_3 + a_5^2 b_0^2\, , \cr
\RR_y  \quad & = \quad
 a_2^2 b_5^2 - a_2 a_4 b_4 b_5
- 2 a_2 a_5 b_2 b_5
+ a_2 a_5 b_4^2
+ a_4^2 b_2 b_5
- a_4 a_5 b_2 b_4
 + a_5^2 b_2^2 \, . \cr} $$
For instance, for $i=3$ and $j=2$ we find
$\,\res^T_f(x^3 y^2) = P_{32}/ \RR_\infty^2$, where
$$ \eqalign{ P_{32} \,\,\,= \,\,\,
&\,\phantom{-} \,      a_0^2 a_1 b_2^2  b_4
\,-\, 2 a_0^2 a_2 b_1 b_2 b_4
\,+\, a_0^2 a_2 b_2^2 b_3
\,-\, a_0^2 a_3 b_2^3
\,+\, a_0^2 a_4 b_1 b_2^2
\,-\, a_0 a_1^2 b_2^2 b_3   \cr &
\,+\, 2 a_0 a_1 a_2 b_1 b_2 b_3
\,-\, 2 a_0 a_1 a_4 b_0 b_2^2
\,+\, 2 a_0 a_2^2 b_0 b_1 b_4
\,-\, 2 a_0 a_2^2 b_0 b_2 b_3   \cr &
\,-\, a_0 a_2^2 b_1^2 b_3
\,+\, 2 a_0 a_2 a_3 b_0 b_2^2
\,+\, a_1^2 a_3 b_0 b_2^2
\,-\, a_1 a_2^2 b_0^2 b_4
\,-\, 2 a_1 a_2 a_3 b_0 b_1 b_2    \cr &
\,+\, 2 a_1 a_2 a_4 b_0^2 b_2
\,+\, a_2^3 b_0^2 b_3
\,-\, a_2^2 a_3 b_0^2 b_2
\,+\, a_2^2 a_3 b_0 b_1^2
\,-\, a_2^2 a_4 b_0^2 b_1 \, .
\cr} \eqno (0.3) $$
It is convenient to review the toric algorithm of
[CD] for computing global residues by means of this example.
First introduce the homogeneous polynomials
$\,F_s(x,y,z)\, :=$ \break $ \,z^2 \cdot f_s({x /z},{y /z})\,$ 
for $s=1,2$.
Next consider the following meromorphic $2$-form on ${\bf P}^2$:
$$ { x^{i-1} y^{j-1} \over z^{i+j-3} F_1\, F_2} \cdot \Omega \,, 
\eqno (0.4) $$
where $\,\Omega \,=\,  x dy \wedge dz -  y dx \wedge dz + 
z dx \wedge dy \,$  denotes the Euler form on ${\bf P}^2$.
The residue (0.2) in the torus $({\bf C}^*)^2$ coincides with the
toric residue of (0.4) in ${\bf P}^2$.

Suppose, for simplicity, that $i\geq 1$, $j\geq 1$ and 
$i+j > 3$.  Consider the
homogeneous ideal $I = \,\langle z^{i+j-3}, F_1, F_2 \rangle$ in
the polynomial ring $K[x,y,z]$ over the field
$K = {\bf Q}(a_0,a_1,\ldots,b_5)$.  
The quotient modulo this ideal is a one-dimensional  
$K$-vector space
in the socle degree $i+j-2$. The homogenized Jacobian
$\,J(x,y,z) := (z^{i+j}/xy)\,J^T(x/z,y/z)$ has degree $i+j-2$ and is
non-zero modulo $I$.  Thus, the monomial $x^{i-1}y^{j-1}$ may 
be written
as $\lambda\,J(x,y,z)$ modulo $I$, where $\lambda \in K$.  
The desired residue 
$\res^T(x^iy^j)$ is then given by $4  \lambda$.  
The coefficient $\lambda$
may be computed, for example, as the ratio
 of the normal form
of $x^{i-1}y^{j-1}$ and the normal form of
$J$ relative to a Gr\"obner basis of $I$.

To prove a denominator formula like (0.2') we use the
following technique. We replace the form
$z^{i+j-3}$ by a generic homogeneous polynomial $F_0(x,y,z)$
of degree $i+j-3$. Note that
$F_0$ has ${i+j-1 \choose 2}$ indeterminate
coefficients, say, $c_0,c_1,c_2,\ldots$.
Consider  the $2$-form
$$\, { x^{i-1}y^{j-1} \over F_0 F_1 F_2} \cdot \Omega \,. 
\eqno (0.5) $$
Now all three forms in the denominator of (0.5)
are generic relative to their degrees.
In~\S~1  we study this situation for an arbitrary projective toric
variety in the role of ${\bf P}^2$. Theorem~1.4 implies that
the denominator of the toric residue of (0.5) equals
the resultant ${\cal R} = {\cal R}(F_0,F_1,F_2)$.
We now apply the specialization
$\,F_0 \mapsto z^{i+j-3} $, which sets all
 but one of the variables $c_0,c_1,\ldots$ to zero.
It takes  (0.5) to  (0.4), and by Lemma~3.4,
it takes $\,{\cal R}\,$ to $\,R_\infty^{i+j-3}$, as desired.
Such a specialization from a generic polynomial $F_0$
to a monomial will connect residues in the torus (\S 3) to
toric residues (\S 1). This technique
will reduce Theorem~3.2 to Theorem~1.4.

In \S 2 we express the sparse resultant as the determinant
of a Koszul-type complex which involves the Jacobian. In some
special cases (Corollary~2.4) we obtain Sylvester-type formulas
which generalize the approach in [GKZ, \S III.4.D] (see also [Ch]).

\medskip
\noi {\bf Acknowledgements:\ } We are  grateful to David Cox,
Fernando Cukierman,
Irena Peeva, and Richard Stanley for their very helpful suggestions.
Part of the work on this paper was done while Eduardo Cattani
was visiting the University of Grenoble and the University of Buenos
Aires; he is thankful for their support and hospitality.
\eject
\noi {\bf \S 1. Residues, Jacobians and resultants 
in toric varieties}
\smallskip
\noindent
We begin with a review of basic concepts from 
toric geometry including the {\sl toric residue}.  
For details and proofs see
[F],[O],[C1],[C2], and [CCD].
Let $X = X_P$ denote the projective toric variety 
defined by an integral, $n$-dimensional polytope
 $$P \,\,:=\,\, \bigl\{ \,m 
\in \R^n \ : \ \ip m {\eta_i} \geq - \bi \,\,
\, \hbox{for} \,\,\, i=1,\ldots,s \, \bigr\}\, , \eqno{(1.1)}$$
where the $\eta_i$ are the first integral vectors in the inner
normals to the facets of $P$. Thus, $X$ is the toric variety 
associated 
with the lattice $M=\Z^n$ and the inner normal fan $\Sigma(P)$
as in [F, \S 1.5].
 We introduce the polynomial ring  $S := \C[x_1,\ldots,x_{s}]$,
where the variable $x_i$ is associated to the
generator $\eta_i$ and hence to a torus-invariant 
irreducible divisor
$D_i$ of $X$.  The Chow group $A_{n-1}(X)$ of
invariant Weil divisors is presented by the exact sequence
$$0 \to  M  \to \Z^{s} \to A_{n-1}(X) \to 0\, \eqno (1.2) $$
where the left morphism sends $m\in M$ to the $s$-tuple
$\,\ip m {\eta}\,\,:= \,\,(\ip m {\eta_1},\ldots,\ip m {\eta_{s}})$.

Let $Z $ denote  the algebraic subset of
$\C^{s}$ defined by the radical monomial ideal 
$$ \langle \prod_{\eta_i \not\in \sigma} x_i \ , \ \sigma \hbox
{ a cone of } \Sigma(P) \ \rangle
\quad \subset \quad S  \, . $$
  The algebraic group 
$G:={\rm Hom}_{\Z} ( A_{n-1}(X), \C^*) \hookrightarrow (\C^*)^s$
acts naturally on $\C^s$ leaving $Z $ invariant.  
The toric variety
$X$ may be realized as the categorical quotient of 
$\C^s \backslash Z $ by $G$
(see [C1]).  When $X$ is simplicial (i.e.~$P$ is simple), 
then the $G$-orbits are closed and $X$ is
the geometric quotient of $\C^s \backslash Z $ by $G$.  
The torus $(\C^*)^s $ lies in $ \C^s \backslash Z $ 
and maps onto the dense torus in $X$ under the quotient map.
    
Given $a\in \N^s$ we write $x^a$ for the monomial 
$\Pi_{i=1}^{s} x_i^{a_i}$.  
As in [C1] the right morphism in (1.2)
defines an $A_{n-1}(X)$-valued  grading of 
the polynomial ring $S$:
$$ \deg(x^a) \,\,:= \,\,[\sum_{i=1}^{s} a_i\,D_i]
\,\,\,\in \,\,\, A_{n-1}(X) \, . \eqno (1.3) $$ 
Let $S_\alpha$ denote the graded component
of $S$ of degree $\alpha$.  We abbreviate
$\beta_0 := [\sum_i D_i]$ and
$\beta := [\sum_i \bi D_i]\in A_{n-1}(X)$.  The divisor $\beta$ is
ample  and $\,S_\beta \cong H^0(X,\LL)$,
where $\LL=\ox(\beta)$ is
the line bundle associated to $\beta$ (see [F, \S 3.4]).
Thus, a homogeneous polynomial $F$ of degree $k\beta$ represents
a global section of $\LL^{k}$, and we may
consider its zero set in $X$.

A monomial $x^a$ has degree $k\beta$, $k\in \N$, if and only if
there  exists $m(a)\in \Z^n$ such that 
$$\ip {m(a)} {\eta_i} + k\bi \,\,\, = \,\, \, a_i\quad \,
\hbox{for } \quad i=1,\ldots,s\,.$$
The point $m(a)$ is unique and, since
$a_i\geq 0$, it lies in $ kP\cap\Z^n$.  Therefore, the map
$$ kP\cap\Z^n \, \to \, S_{k \beta} \, ,\quad
m \,\mapsto \,\prod_{i=1}^s \,x_i^{\ip {m} {\eta_i} +
k\bi}\eqno{(1.4)}$$
defines a bijection between integral points in $kP$ and
monomials of degree $k\beta$  
or, equivalently, between Laurent polynomials supported in $kP$ and 
homogeneous polynomials of degree $k\beta$ in $S$. 
If $f(t_1,\ldots,t_n)$ is supported in $kP$ then its 
image is the $kP$-{\sl homogenization}
$$F(x_1,\ldots,x_s) \,\, \, =  \,\,\,
\bigl(\prod_{i=1}^s \,x_i^{k\bi}\bigr)\cdot
f(t_1(x),\ldots,t_n(x))  \,\, \in \,\, S_{k \beta}\, , \eqno(1.5)$$
$$\hbox{ where} \qquad
 t_j(x) \quad = \quad \prod_{i=1}^s \,x_i^{\ip {e_j} {\eta_i} }\quad \qquad
(j=1,\ldots,n) $$
and  $\{e_1,\ldots,e_n\}$ is the standard basis of $\Z^n$.
By restricting (1.4) we also get a bijection between monomials
$x^a$ of degree $k\beta - \beta_0$ and integral points in
$(kP)^\circ$, the interior of $ k P $.

\bigskip

\proclaim Proposition~1.1.    The ring 
$S_{*\beta}= \bigoplus_{k=0}^\infty S_{k\beta}$ is
Cohen-Macaulay of dimension $n+1$, with canonical module 
$\,\omega_{S_{*\beta}} = \bigoplus_{k=0}^\infty S_{k\beta-\beta_0}$.
Fix positive integers $k_0,\ldots , k_n$ and let
$\kappa = k_0+\cdots+k_n$, $\rho = 
\kappa \beta - \beta_0$. Given
$F_i\in S_{k_i\beta}$ for $i=0,\ldots,n$ such that $F_0, 
\ldots, F_n$ 
have no common zeroes in $X$, then:
\item{(i)} $F_0,\dots,F_n$ are a regular sequence in $S_{*\beta}$ 
and,
hence, in $\omega_{S_{*\beta}}$.
\smallskip
\item{(ii)} The degree $\rho$ component  $R_\rho$ of the quotient
$R = S_{*\beta}/\langle F_0,\dots, F_n\rangle$ has 
${\bf C}$-dimension~$1$.  

\noi{\bf Proof: } See [B, Theorem~2.10 and
Proposition~9.4] and [C2, Proposition~3.2].\ \ $\diamond$

\bigskip

We next recall the construction  of the
Euler form $\Omega$ and the toric Jacobian $J(F)$ 
(see [BC,\S9], [C2,\S4]).
For any subset 
$I = \{i_1,\dots, i_n\}$ of $\{ 1,\ldots, s\}$ we abbreviate
$$ \det(\n_I) \,\,:=\,\, \det(\langle e_\ell,\n_{i_j}
        \rangle_{\scriptscriptstyle{1 \le \ell,j \le n}})\,, \ \ 
dx_I = dx_{i_1}\wedge\cdots\wedge dx_{i_n}\,, \ \ 
\hat{x}_I = \Pi_{j \notin I} x_j \,. $$  
Note that  the product
$\det(\n_I) dx_I$ is independent of the ordering of
$ i_1,\dots, i_n$.  The {\it Euler form} on $X$ is the following
sum over all $n$-element subsets $I \subset \{ 1,\ldots, s\}$:
$$\Omega \quad := \quad 
\sum_{|I| = n} \det(\n_I)\, \hat{x}_I\, dx_I \, . $$
The Euler form $\Omega$ may be 
characterized by the property that $\Omega/(x_1\cdots x_s)$ is 
the rational extension to $X$ of the $T$-invariant holomorphic form 
$\frac {dt_1} {t_1} \wedge \cdots \wedge \frac {dt_n} {t_n}\,$
on the torus $T$.

\medskip

As in Proposition~1.1, 
consider homogeneous polynomials $F_0,F_1,\ldots,F_n$ 
where  $\deg(F_i)$ $ = k_i \beta $
and $\kappa = k_0+\cdots+k_n$. Then there exists
a polynomial $J(F) \in S_{\kappa \beta - \beta_0}$ such that
$$\textstyle{\sum_{i=0}^n} (-1)^i F_i \cdot  dF_0\wedge \cdots 
\wedge
dF_{i-1} \wedge dF_{i+1} \wedge \cdots \wedge dF_n \quad = 
\quad J(F)\cdot\Omega\, .\eqno(1.6)$$
Furthermore, if $I = \{ i_1,\dots, i_n\}$ is such that
$\eta_{i_1},\dots,\eta_{i_n}$ are linearly independent, then
$$ J(F)  \quad = \quad \frac 1 {\det(\eta_I)\,\hat x_{I}} \ 
\det\pmatrix{k_0F_0&k_1F_1&\ldots&k_nF_n\cr
{\partial F_0/\partial x_{i_1}} &{\partial F_1/\partial x_{i_1}}
&\ldots& 
{\partial F_n/\partial x_{i_1}}\cr
\vdots & \vdots & \ddots & \vdots \cr
{\partial F_0/\partial x_{i_n}} &{\partial F_1/\partial x_{i_n}}
&\ldots& 
{\partial F_n/\partial x_{i_n}} \cr
}\, . \eqno(1.7)$$
The polynomial $J(F)$ is called the {\it toric Jacobian}
of $F = (F_0,F_1,\ldots,F_n)$.
\medskip

In the special case $k_0 = k_1 = \cdots = k_n = 1$ 
the toric Jacobian
can also be computed as follows.  Let $f_0,\ldots,f_n$ be Laurent
polynomials supported in $P$ and let $F_0,\ldots,F_n$ denote their
$P$-homogenizations as in (1.5).  
Let $P\cap\Z^n = \{m_1,\ldots,m_\mu\}$ and 
$$j(t) \quad :=  \quad
\det \pmatrix{
f_0 &
f_1 & \dots &
f_n  \cr
t_1 {\partial f_0 \over \partial t_1} &
t_1 {\partial f_1 \over \partial t_1} & \dots &
t_1 {\partial f_n \over \partial t_1}  \cr
\vdots & \vdots & \ddots & \vdots \cr
t_n {\partial f_0 \over \partial t_n} &
t_n {\partial f_1 \over \partial t_n} & \dots &
t_n {\partial f_n \over \partial t_n}  \cr } \, . \eqno (1.8) $$

\bigskip

\proclaim Proposition~1.2.    Let $f_j = \sum_{i=1}^\mu\
u_{ji}\,t^{m_i}$ and set $\tilde m_i = (1,m_i)\in \Z^{n+1}$. Then,
$$j(t) \quad = \,\,\,
\sum_{1 \leq i_0 < i_1 < \ldots < i_n \leq \mu} \!\! \!\!
[i_0 i_1 \ldots i_n ]
\cdot \det( {\tilde m}_{i_0},{\tilde m}_{i_1},\ldots,
{\tilde m}_{i_n}) 
\cdot
{t}^{ { m}_{i_0} + { m}_{i_1} + \cdots +{ m}_{i_n}} \,,$$
where the brackets denote the maximal minors of the 
coefficient matrix:
$$ [i_0 i_1 \ldots i_n ] \quad := \quad
\det \pmatrix{
 u_{0 i_0} &  u_{0 i_1} & \dots &  u_{0 i_n} \cr
   \vdots  &    \vdots  & \ddots &    \vdots  \cr
 u_{n i_0} &  u_{n i_1} & \dots & u_{n i_n} \cr } \, .  $$
Moreover, $j(t)$ is supported in $((n+1)P)^\circ$ and its
$(n+1)P$-homogenization is $ x_1\cdots x_s J(F)$.

\noindent {\bf Proof:}  We consider the $(n+1)\times\mu$ matrix
$\tilde A = (\tilde m_1,\ldots,\tilde m_\mu)$, the
 $\mu\times\mu$ diagonal matrix 
$D={\rm diag}(t^{m_1},\ldots,t^{m_\mu})$ and
the  $\mu\times (n+1)$ matrix $U$, obtained by 
transposing the matrix
of coefficients $(u_{ji})$.  Their product  
$\tilde A\cdot D \cdot U$
equals the $(n+1)\times (n+1)$ matrix in (1.8).  
The first assertion amounts
to the {\sl Cauchy-Binet formula\/} for 
$\, j(t) \, = \, (\wedge_{n+1} \tilde A ) \cdot
(\wedge_{n+1} D ) \cdot
(\wedge_{n+1} U ) $.
\smallskip
\noindent
If the sum ${ m}_{i_0} + { m}_{i_1} + \cdots +{ m}_{i_n}$ lies in the
boundary of $(n+1)P$, then all $m_{k_j}$ lie in a facet of $P$ and
the determinant 
$ \det( {\tilde m}_{i_0},{\tilde m}_{i_1},\ldots,{\tilde
m}_{i_n}) $ must vanish.  Consequently, $j(t)$ is supported in the
interior of $(n+1)P$.
The final statement follows from (1.6) together with 
$$j(t)\ \frac {dt_1} {t_1} \wedge \cdots \wedge \frac {dt_n} {t_n} 
\quad = \quad
{\sum_{j=0}^n} (-1)^j f_j\, df_0\wedge \cdots \wedge
d f_{j-1} \wedge d f_{j+1} \wedge \cdots \wedge df_n\, .
\ \ \diamond
$$

\bigskip 
We now return to general $k_0,\ldots,k_n$. Suppose
that $\Efes$  have no common zeroes in $X$. Then
$\,R_\rho \cong \C\,$ by (ii) in Proposition~1.1.
In [C2] Cox constructs an explicit isomorphism
$\,\res^X_F \colon R_\rho \rightarrow \C $
whose value on the
toric Jacobian is the positive integer
$$ \res^X_F(J(F)) \quad = \quad  \left( \prod_{j=0}^n k_j \right)
\cdot n! \cdot {\rm vol}(P) \, ,  \eqno{(1.9)} $$
where vol$(\, \cdot \, )$ denotes the standard volume in 
${\bf R}^n$.
The isomorphism $\,\res^X_F( \,\cdot\, ) \,$ is called the 
{\sl toric residue}.
From (1.9) we conclude that
$$ J(F) \, \hbox{ defines a non-zero element in } R_\rho \,. 
\eqno{(1.10)} $$
\medskip
\noi We next present an affine interpretation of the toric residue.
Let $f_j$ be a generic Laurent polynomial with
Newton polytope $ k_jP$.  Let $F_j \in S_{k_j\beta}$ be the
$k_jP$-homogenization of $f_j$.  Given a
homogeneous polynomial $H$ of {\sl critical degree\/} 
$\rho = \kappa\beta - \beta_0$, the expression
$$\frac {H\ \Omega} {F_0\cdots F_n}$$
defines a meromorphic $n$-form on $X$.
Its restriction to $T$ may
be written as
$$\frac {h} {f_0\cdots f_n}\ 
\frac {dt_1} {t_1} \wedge \cdots \wedge \frac {dt_n} {t_n} \, ,$$ 
where $h$ is a Laurent polynomial supported in $(\kappa P)^\circ$.
Our generic choice of $f_0,\ldots,f_n$ guarantees 
(cf.~[K1,\S 2]) the following properties
for each $i=0,\ldots,n$: The finite set 
$\,V_i :=\{x\in X\,:\, F_j(x)
 = 0\,;\,j\not= i\}$ lies in the torus $T$, hence
$V_i=\{t\in T\,:\, f_j(t) = 0\,;\,j\not= i\}$, and
 the function $h/f_i$ is regular at the points of $V_i$.  

The following result is a consequence of Theorem~0.4 in [CCD]:

\proclaim Proposition~1.3. 
For any fixed $i \in \{0,\ldots,n\}$,
the toric residue  equals
$$\res_F^X(H)\quad =\quad (-1)^i\, \sum_{\xi\in V_i} 
\res_\xi \bigl(\frac {h/f_i} {f_0\cdots f_{i-1} f_{i+1} \cdots f_n}\ 
\frac {dt_1} {t_1} \wedge \cdots \wedge \frac {dt_n} {t_n}\bigr) \, .
\eqno{(1.11)}$$
Here the right-hand side is a sum of Grothendieck residues 
{\rm([GH], [T]}; see also \S3{\rm)}
relative to the divisors $\{f_j(t)=0\}\subset T$, $j\not=i$.

\medskip
\noindent {\bf Remarks.}
\smallskip\noi
{\bf i)} Even though Theorem~0.4 in [CCD] is only
stated for simplicial toric varieties, it is valid for 
arbitrary complete toric varieties provided $V_i$ lies in $ T$,
by passing to a desingularization.
\smallskip\noi
{\bf ii)} Note that while the right side of (1.11)
makes sense for every Laurent polynomial $h$,
Proposition~1.3
asserts  that, if  $h$ is supported 
in $(\kappa P)^\circ$, then that expression  is
independent of $i$. 
\bigskip

We next consider $n+1$ polynomials having indeterminate coefficients:
$$ F_i(u;x) \quad := 
\quad \sum_{a\in \AA_{k_i\beta}} u_{ia}\,x^a \qquad
\hbox{for} \quad i=0,\ldots,n \, ,\eqno{(1.12)} $$
where $\,\AA_{k_i \beta} \, := \, \{\, a \in \N^s \, : \,
\deg(x^a) = k_i \beta \, \}  $.
We shall work in the polynomial ring 
$$ C \, := \, A[x_1,\ldots,x_s]  \qquad
\hbox{over } \quad A \, := 
\, \Q\,[\,u_{ia}\,\,;\ i=0,\ldots,n\,;\ a\in \AA_{k_i\beta}] \,  .$$
We endow the polynomial ring $C$ with 
the $A_{n-1}(X)$-grading given by (1.3).  
For any $H \in C_\rho$, the expression (1.11) 
depends rationally on the
coefficients of $F_0,\ldots,F_n$ and hence 
defines an element in the field
of fractions of $A$, which we also denote $\res_F^X(H)$.

As in [GKZ, 3.3; 8.1] we define the {\sl resultant\/} 
associated with the  bundles  $\LL^{k_0}\! ,\ldots,\LL^{k_n}$.  
It is an irreducible polynomial 
$\RR_{\LL^{k_0}\! ,\ldots,\LL^{k_n}}(u)\in
A$  with integral coefficients, uniquely defined up to  sign, 
which vanishes for some specialization of the coefficients if and 
only if the corresponding sections 
$F_0,\ldots,F_n$ have a common zero in $X$.  Via the 
correspondence (1.4) between homogeneous polynomials of degree
$k\beta$ and Laurent polynomials supported in $kP$, the resultant
$\RR_{\LL^{k_0}\! ,\ldots,\LL^{k_n}}(u)$  agrees with the
{\sl mixed sparse resultant\/} (see [PSt],[S2])
associated with the support sets 
$\, k_0P\cap\Z^n,\ldots,k_nP\cap\Z^n$.

The degree of the resultant is computed as follows. Suppose
$ k_0 \geq \ldots \geq k_n$.
Consider the lattice affinely generated
by the integral points in $k_0 P$. It has finite index in $\Z^n:$
$$\ell \quad := \quad \left[\, \Z^n \, : \, 
{\rm aff}_{\Z} ( k_0 P \cap
\Z^n) \, \right] \, . \eqno (1.13)$$
Note that $\ell = 1$ if ${\cal L}^{k_0}$ is very ample. The degree
of $\RR_{\LL^{k_0}\! ,\ldots,\LL^{k_n}}(u)$ in the coefficients
of the $i$-th form $F_i$ equals, by [PSt, Corollary~1.4],
$$k_0 \cdots k_{i-1} k_{i+1} \cdots k_n 
\cdot n \, ! \cdot { 1 \over \ell} \cdot \vol(P) \, .  \eqno (1.14)  $$

\smallskip

\noi We now state and prove the main result of this section:

\medskip

\proclaim Theorem~1.4. For any $H\in C_\rho$,
the product 
$ \, \RR_{\LL^{k_0}\! ,\ldots,\LL^{k_n}}(u) \cdot \res^X_F(H)
\,$ lies in $\, A$.

\noi
{\bf Proof:} As noted above, for values of $u$ in a Zariski open set,
$\Efes$ have no common zeroes in $X$ and, for every 
$i=0,\ldots,n$, the set
$V_i = \{x\in X\,:\,F_j(x)=0,\, j\not=i\}$ is
finite and contained in $T$.  Thus, setting for simplicity $i=0$, we
have, as in (1.11):
$$\res_F^X(H)\quad =\quad \sum_{\xi\in V_0} 
\res_\xi \bigl(\frac {h/f_0} {f_1\cdots  f_n}\ 
\frac {dt_1} {t_1} \wedge \cdots \wedge \frac {dt_n} {t_n}\bigr) \, .
\eqno{(1.15)}$$
We may further assume
that the zeroes of $f_1,\ldots,f_n$ are simple and, therefore,
each term in the right hand side of (1.15) may be written as (see
[GH, page 650]):
$$\res_\xi \bigl(\frac {h/f_0} {f_1 \cdots f_n}\ 
\frac {dt_1} {t_1} \wedge \cdots \wedge 
\frac {dt_n} {t_n}\bigr)
\,\, = \,\, \frac {h(\xi)} {f_0(\xi)\,
\cdot\,J^T_{f_1,\ldots,f_n}(\xi)}
\,\, = \,\,
\frac {a_\xi(u_1,\ldots,u_n)} {f_0(\xi)\,
\cdot\,b_\xi(u_1,\ldots,u_n)}\,,
\eqno (1.16) $$
where $\,J^T_{f_1,\ldots,f_n} = 
\displaystyle{\det(t_j\,\frac {\partial f_i} {\partial t_j})}$, 
the symbol $u_i$ 
stands for the vector $\left( u_{ia} \, : 
\, a \in {\cal A}_{k_i \beta}\right)$
of coefficients of $\, f_i\,$, 
and $\,a_\xi$, $\,b_\xi\,$ are algebraic 
functions in these coefficients.

We now sum (1.16) over all points $\xi$ in $V_0$. 
To get the best possible
denominator even if $\ell > 1$, we must organize the sum (1.15) as
follows. First, we may assume that $P$ contains the origin. Then the
affine lattice agrees with the linear lattice,
$${\rm aff}_\Z (k_0 P \cap \Z^n) \quad = 
\quad {\rm lin}_\Z ( k_0 P
\cap \Z^n)\, , \eqno (1.17)$$
and the inclusion of (1.17) in $\Z^n$ defines a morphism of tori
$\pi: \, T \to \left( \C^* \right)^n$. The map $\pi$ is a finite
cover of degree $\ell$, and the Laurent polynomial $f_0$ is constant
along the fibers of $\pi$. Hence, if 
$\eta= \pi(\xi)$ for $\xi \in V_0$,
then we can define $f_0(\eta) := f_0(\xi)$. Therefore,
$$ \res_F^X(H) \quad = \quad \sum_{\eta \in \pi(V_0)} 
\frac 1 {f_0(\eta)}\,
 \sum_{\xi \in \pi^{-1}(\eta)} 
\frac {a_\xi(u_1,\ldots,u_n)} {b_\xi(u_1,\ldots,u_n)}\,.$$
This expression depends rationally on $u_0, u_1, \ldots, u_n$. 
This implies
$$\res_F^X(H) \quad
=\quad
 \frac {A(u_0,u_1,\ldots,u_n)} {(\prod_{\eta \in \pi(V_0)}
f_0(\eta))\cdot B(u_1,\ldots,u_n)}\, ,$$
where $A$ and $B$ are polynomials.
It follows from [PSt, Theorem~1.1] that
$$ \prod_{\eta \in \pi(V_0)}f_0(\eta) \quad =
\quad \RR_{\LL^{k_0}\! ,\ldots,\LL^{k_n}}(u_0,u_1,\ldots,u_n)
\cdot C (u_1,\ldots,u_n)$$
for some rational function $C$.
Therefore, there exist polynomials $A_0, B_0$ such that
$$\res_F^X(H) \quad
=\quad
 \frac {A_0(u_0,u_1,\ldots,u_n)} {\RR_{\LL^{k_0}\! ,\ldots,
\LL^{k_n}}(u_0,u_1,\ldots,u_n) \cdot B_0(u_1,\ldots,u_n)}\,.$$
Replacing the role played by the index $0$ by any other index
$i=1,\ldots,n$, we deduce that 
$$\res_F^X(H) \quad
=\quad
 \frac {P(u_0,u_1,\ldots, u_n)} {\RR_{\LL^{k_0}\! ,\ldots,
\LL^{k_n}}(
u_0,u_1,\ldots,u_n)}$$
for some polynomial $\, P \in A$. \ \ $\diamond $

\bigskip
\noi
{\bf Remark 1.5.}  Suppose $P$ is the standard simplex in $\R^n$.  
Then
$X\cong \P^n$, $\beta$ is the hyperplane class,
$s=n+1$, and
$F_j(x_0,\ldots,x_n)$ is a homogeneous polynomial of degree $k_j$. 
The assumption that $\Efes$ have no common zeroes in $\P^n$ means
that their only common zero in $\C^{n+1}$ is $0$.
For any homogeneous polynomial $H$ of degree $\rho=\kappa - (n+1)$,
the toric residue $\res_F^{\P^n}(H)$ associated with the 
$n$-rational form
$\frac {H} {F_0\cdots F_n} \Omega $ on $\P^n$, coincides 
([PS], [CCD,\S5]) with
the Grothendieck residue at the origin of $\C^{n+1}$ of the 
$(n+1)$-form
$$\frac {H} {F_0\cdots F_n} \ dx_0\wedge dx_1 \wedge 
\cdots \wedge dx_n \, .$$
In this situation, it has been observed by 
Ang\'eniol [A] that Theorem~1.4
follows from the work of Jouanolou (see, for example, [J1, 3.5]).

\vskip 1cm
\eject

\noi {\bf \S 2. Jacobian formulas for the sparse resultant}
\smallskip
\noi
Let $F_0,\ldots,F_n$  be generic forms as in (1.12), let
$A$ be the polynomial ring on their coefficients,
and let $C = A[x_1,\ldots,x_s]$ be graded 
by the Chow group $A_{n-1}(X)$ via (1.3). 
The given forms together with their
toric Jacobian $J(F)$ define a map of free  $A$-modules
$$ \eqalign{
\Phi \, : \,\, C_{\rho-k_0\beta}\times\cdots\times 
C_{\rho-k_n\beta}\times A
\quad & \rightarrow \quad \,\, C_{\rho}\, , \cr
(\,\Lambda_0\, ,\, \ldots \, ,\, \Lambda_n,\,\Theta\, ) 
\quad \, & \mapsto 
\,\,\,   \sum_{i=0}^n \Lambda_i\, F_i + \Theta\,J(F) \,.\cr} 
\eqno (2.1) $$
For any particular choice of complex coefficients $\,u = c \,$
we abbreviate $\, F_i^c(x) := F_i(\,c \,; \,x \,)$.
The resultant $\,{\cal R} = {\cal R}_{{\cal L}^{k_0}, \ldots, 
{\cal L}^{k_n}}
\in A\,$ considered in Theorem~1.4 satisfies
$\,{\cal R}(c) = 0 \,$ if and only if the forms
$F_0^c,\ldots,F_n^c$ have a common zero in the toric variety $X$. 
Let
$$\Phi_{c} \,\, : \,\, 
S_{\rho-k_0\beta}\times\cdots\times S_{\rho-k_n\beta}
\times \C \,\rightarrow \, S_{\rho}\eqno{(2.2)}$$
denote the ${\bf C}$-linear map derived from (2.1)
by substituting $c$ for $u$.

\medskip

\proclaim Proposition~2.1. The map $\Phi_{c}$ is surjective if and 
only if $\RR(c) \not=0$.

\noindent {\bf Proof:}  For the if direction suppose
$\RR(c)\not =0$. Then $\,F_0^c,\ldots,F_n^c\,$ have
no common zeroes in $X$.
Proposition~1.1 (ii)  together with (1.10) implies
the surjectivity of $\Phi_{c}$.

For the converse, let ${\cal V}$ denote the affine variety 
in the space of coefficients consisting of all $\,c \,$ such
that the polynomials $\,F_0^c,\ldots,F_n^c\,$ have a common zero in
the torus $(\C^*)^s$.  
Fix $c\in {\cal V}$ and let  $p\in(\C^*)^s$ be such a common zero.  
It follows from (1.7) that
$ x_1 x_2 \cdots x_s\cdot J(F)$ lies in the ideal generated by 
$F_0,\ldots,F_n$ in $S$ and hence $J(F)$ vanishes at $p$.  
If a monomial
$x^a$ of degree $\rho$ were in the image of $\Phi_{c}$
then  $x^a(p) =0$ which is impossible.  Thus, for $c\in {\cal
V}$, $\Phi_{c}$ is not surjective. We conclude that ${\cal V}$ is 
contained in the algebraic variety defined by the vanishing of all
maximal minors of  $\Phi_{c}$.  Since the closure of ${\cal V}$ is
the locus where the resultant ${\cal R}$ vanishes, 
the only if-direction follows.
 \ \ $\diamond$

\bigskip
\noi For any subset $J \subseteq \{0,\ldots, n\}$ we set
$\, k_J := \sum_{i \in J} k_i$.
For $\,0\leq j \leq n+1 \,$ denote 
$$W_j \quad := \quad \bigoplus_{|J| = j} C_{k_J \beta - \beta_0} \, .
\eqno{(2.3)}$$
From the  Koszul complex on $\Efes$ we derive the
following complex of free $A$-modules:
$$ 0 \, \longrightarrow W_0 \, \mapname{\varphi_0} \, W_1  
\mapname{\varphi_1}
\, \ldots \, \mapname{\varphi_{n-1}} \, W_n \, \mapname{\varphi_n}
\, W_{n+1}  \, \longrightarrow \, 0 \, . \eqno{(2.4)}$$
This construction is an instance of [GKZ, \S 3.4.A].
Note that $\,W_0 = 0$, $\,W_{n+1} = C_\rho$, and 
$\,W_n =  C_{\rho-k_0\beta}\times\cdots\times C_{\rho-k_n\beta}\,$.
Define $(\varphi_{n-1},0) : W_{n-1} \,
\longrightarrow W_n \oplus A$  by adding $0$ in the coordinate
corresponding to $A$, and consider the modified complex
$$ 0 \,\longrightarrow W_1 \, \mapname{\varphi_1} \, W_2  
\mapname{\varphi_2}\,
\, \cdots \,\,  \mapname{\varphi_{n-2}}\, W_{n-1} 
 \mapname{(\varphi_{n-1},0)} \, W_n\oplus A \, \mapname{\Phi}
\, W_{n+1}  \, \longrightarrow \, 0 \, . \eqno{(2.5)}$$
For  any particular choice of coefficients $u=c$ 
in (2.5) we get a complex of ${\bf C}$-vector spaces:
$$ 0 \longrightarrow \bigoplus_{i} S_{k_i \beta - \beta_0} 
\mapname{\varphi_1^c}  \, \cdots \, \mapname{(\varphi_{n-1}^c,0)} 
\bigoplus_{|J|=n} S_{k_J \beta -\beta_0}\times \C
\mapname{\Phi_c} S_{\rho} \longrightarrow 0 \, .
\eqno{(2.6)} $$ 
Let $D$ denote the determinant (see [GKZ, Appendix A]) 
of the complex 
of $A$-modules (2.5) with respect a fixed choice of monomial bases
for the $A$-modules $\, W_1,\ldots, W_{n+1}$. This is an element in
the field of fractions of $A$.  We shall prove that 
it is a polynomial
in $A$. Suppose $ k_0 \geq \ldots \geq k_n$ and let $\ell$ be the
lattice index defined in (1.13).

\proclaim Theorem~2.2.
\item{(i)} The complex of $\C$-vector spaces (2.6) is exact if 
and only 
if $\RR(c) \not= 0$.
\item{(ii)} The determinant $D$ of the complex (2.5) equals the
greatest common divisor of all (not identically zero)
maximal minors of a matrix representing the
$A$-module map  $\Phi$. 
\item{(iii)}The determinant  $D$ equals  $\RR^\ell$. 
\item{(iv)} If ${\cal L}^{k_0}$ is very ample then the resultant 
$\RR$ 
may be computed as the greatest common divisor of all 
maximal minors of
any matrix representing $\Phi$.  

\medskip

\noi{\bf Proof: }
We first prove the if-direction in part (i).
Let $\beta$ be an ample divisor and 
 $F_0^c,\ldots,F_n^c$  homogeneous polynomials 
of respective degrees
$k_i  \beta$ without common zeroes in $X$, i.e. 
such that $\RR(c)\not=0$.
By Proposition~1.1 (i),  $F_0^c,\ldots,F_n^c$ is a regular sequence 
in $S_{*\beta}$ and in $\omega_{S_{*\beta}}$ ; consequently, 
the corresponding Koszul complex is acyclic 
[BH, page 49].  Setting $I = \langle F_0^c,\ldots,F_n^c\rangle$ this implies
that
$$ 0 \longrightarrow 
\bigoplus_{i} S_{k_i \beta - \beta_0} \mapname{\varphi_1^c} 
\, \cdots \, \mapname{\varphi_{n-1}^c} 
\bigoplus_{|J|=n} S_{k_J \beta -\beta_0} 
\mapname{\varphi_n^c} S_{\rho} 
\longrightarrow  S_{\rho} / I_{\rho} \longrightarrow 0 
\eqno (2.7)$$
is an exact sequence of ${\bf C}$-vector spaces.
Proposition~2.1 implies that $\Phi_{c}$ is surjective.  
Also, by (1.10),
$\Phi_c( \lambda_1,\ldots,\lambda_n, \theta) = 
\sum_i \lambda_i F_i + \theta J(F) = 0$ implies $\theta =0$. 
These two facts imply that (2.6) is exact.
For the converse of (i) suppose $\RR(c) =0$. Then the map $\Phi_c$ 
is not surjective by Proposition~2.1, and hence (2.6) is not exact.

We next prove part (ii). We claim that $F_0,\ldots,F_n$ is a 
homogeneous regular sequence in the graded Cohen-Macaulay
ring $\, C_{*\beta}:=  \bigoplus_{k=0}^\infty C_{k\beta}$.
We extend scalars and consider $C_{*\beta} \otimes_{\Q} \C$ instead.
Let $N$ be the total number of terms in $\Efes$. The spectrum of
$C_{*\beta} \otimes_{\Q} \C$ equals affine space $\C^N$ times
the $(n+1)$-dimensional affine toric variety 
${\cal X}_\beta:= {\rm Spec}(S_{*\beta})$.  
Let $\cal V$ denote the algebraic set defined by $\Efes$ in 
$\C^N \times {\cal X}_\beta$. 

We shall prove that $\cal V$ has codimension 
$n+1$, by describing the two irreducible components of $\cal V$.
Let $O$ be the origin in ${\cal X}_\beta$ and $\cal M$ 
its maximal ideal.
Hence $\cal M$ is spanned by all non-constant monomials 
in $S_{*\beta}$.
For any $ i \in \{0,\ldots, n\}$, the $x$-monomials 
appearing in $F_i$ all lie
in ${\cal M}^{k_i}$, and their radical equals $\cal M$. 
In other words, 
$F_i(p) \not= 0$ for all $p \in {\cal X}_\beta \setminus  \{O\}$.
Consider the projection from $\C^N \times {\cal X}_\beta$ onto its
second factor and let $\pi$ denote its restriction
to $\cal V$. For $p \in {\cal X}_\beta \setminus  \{O\}$, the fiber
$\pi^{-1}(p)$ is a linear subspace of codimension $n+1$
in $\C^N \times \{p\}$. The fiber $\pi^{-1}(O)$ equals 
$\C^N \times O$,
which has codimension $n+1$ in $\C^N \times {\cal X}_\beta$. We
have shown that codim$({\cal V}) = n+1$, as desired. 

Since $C_{*\beta} \otimes_{\Q} \C$ is graded and Cohen-Macaulay,
we may conclude that $F_0,\ldots,F_n$ is a regular sequence.
The Koszul complex on $F_0,\ldots,F_n$ is exact, and therefore
(2.4) and (2.5) are exact sequences of $A$-modules  
except at $W_{n+1}$. By Theorem~34 in [GKZ, Appendix~A],
the determinant $D$ equals the greatest
common divisor of all maximal minors of $\Phi$.

\smallskip

Part (iv) of Theorem~2.2 follows directly from (ii) and (iii)
and the observation that $\ell = 1$ if ${\cal L}^{k_0}$ 
is very ample.
It remains to prove part (iii). Part (i) implies that
$D(c) = 0 $ if and only if $\RR(c)=0$. 
We also deduce from the irreducibility
of the resultant that $D$ is a power of $\RR$. 
In order to prove $\, D = \RR^\ell \, $, we must show that
the total degree of $D$ equals 
$$ \ell \cdot \deg(\RR) \quad = \quad
\bigl(\, \sum_{i=0}^n k_0 \cdots k_{i-1} k_{i+1} \cdots k_n \bigr)
\cdot n \, ! \cdot \vol(P) \, .
\eqno (2.8)  $$
Let us consider the {\sl Erhart polynomial} for the interior of $P$:
$$ p(j) \quad := \quad
|(j P)^\circ \cap \Z^n| \quad = \quad
 \vol (P) \cdot j^n \,+\, \sum_{i =0}^{n-1} a_i j^i\,.$$
The rank of the free $A$-module $W_j$ equals
$\, \sum_{|J| = j} p(k_J)$.
Taking into account the fact that any non-zero 
maximal minor of $\Phi$ has
to involve the last column and $\deg(J(F))=n+1$ 
in the coefficients of 
$F_0,\ldots,F_n$,   we deduce from Theorem~14 in Appendix
A in [GKZ] that
$$ \eqalign{
& \deg(D)   \quad  = \quad
 \sum_{j=0}^{n+1} (-1)^{n+1-j} \cdot j \cdot \left( 
\sum_{|J|=j} p(k_J) \right) \quad = \cr
 & \, \vol (P) \cdot \!
 \underbrace{\left(\sum_{j=0}^{n+1} (-1)^{n+1-j} \cdot j \cdot 
\! \sum_{|J|=j}\! k_J^n \right)}_{\hbox{$\gamma_n$}} 
\,\,+ \,\,\, \sum_{i=0}^{n-1} a_i  \cdot \!
 \underbrace{\left(
\sum_{j=0}^{n+1} (-1)^{n+1-j} \cdot j \cdot  \!\! \sum_{|J|=j}\!
k_J^i  \right)}_{\hbox{$\gamma_i$}} \,. \cr} \eqno (2.9) $$
To prove the equality of (2.8) and (2.9), it suffices to show the
combinatorial identities:
$$ \gamma_n \,\, = \,\,\, n \,! \cdot (\sum_{j=0}^{n+1} 
\prod_{\nu \not= j} k_\nu) \quad \hbox{ and } \quad   \gamma_i \,\, 
= \,\, 0 
\quad \hbox{for} \,\,\, 0 \leq i \leq n-1\,. \eqno (2.10) $$
Following a suggestion made to us by Richard Stanley, 
we prove a more
general identity:

\proclaim Lemma~2.3. Let $u_{i,j}$ be indeterminates  indexed by
 $i=0,\ldots,n$ and $j=0,\ldots,r$. Then
$$ \sum_{I \subseteq \{0,1,\ldots,n\}} \!\!\!
(-1)^{ |I|} \prod_{j=0}^r \, \bigl( \sum_{i \in I } 
u_{i,j}  \bigr) \qquad
= \qquad \, (-1)^{n+1} \!\!\!\!\!\! \sum_{\phi : \{0,\!\ldots\!,r \}
\rightarrow  \{0,\! \ldots \!,n\} \atop {\rm surjective}}
 \prod_{j=0}^r u_{\phi(j),j} \,. $$

\noindent {\sl Proof:}
The terms in the expansion of the left side correspond to maps
from $\{0,\ldots,r\}$ to subsets $I$ of $\{0,\ldots,n\}$.
Any term which appears at least twice gets cancelled.
What remains are the terms corresponding to surjective maps from
$\{0,\ldots,r\}$ to the full set $I = \{0,\ldots,n\}$.\ \ $\diamond$

\vskip .2cm

\noi We are interested in the special case
 $\,u_{i,0} = 1 \,$ for $ 0 \leq i \leq n$ and
$\,u_{i,j} = k_i \,$ for $ 0 \leq i \leq n$ and $1 \leq j \leq r$.
Under this specialization, Lemma~2.3 implies (2.10) 
and hence part (iii).
This completes the proof of Theorem~2.2.
\ \ $\diamond$

\bigskip
Theorem~2.2 expresses the $\ell$-th power of the resultant 
as an alternating 
product of determinants. Of particular interest are 
those cases when 
one determinant is involved. Such formulas are
called {\it  Sylvester-type}. They have been
studied systematically by Weyman and Zelevinsky [WZ]
in the case when $X$ is a product of projective spaces.

\proclaim Corollary~2.4. Suppose that $(n-1) P$ has no 
interior lattice 
points and either 
\item{(a)} $k_0 = \cdots = k_n = 1, \quad$ or 
\item{(b)} $n P$ has no interior lattice points and 
$k_0 + \cdots + k_n =
n+2, \quad$ or 
\item{(c)} $n=2$ and $P$ is a primitive triangle and
$k_0,k_1,k_2 \leq 2$. 
\item{} Then, the matrix of $\Phi$ is square and $\RR^\ell = 
\det(\Phi)$.

\medskip

Let us discuss the formulas in Corollary~2.4 for the 
case of toric surfaces $(n=2)$. Suppose $k_0 = k_1 = k_2 = 1$ and
the polygon $P$ has no interior lattice points.
Then the matrix of $\Phi$ is square and
$\, \RR\, = \,\det (\Phi)$.
A lattice polygon $P$ has no interior lattice points
if and only if $\,(X,\beta) \,$
is either the Veronese surface in ${\bf P}^5$ or
any rational normal scroll (Hirzebruch surface).
In the former case we recover Sylvester's formula 
for the resultant of three
ternary quadrics [GKZ, \S 3.4.D].
In the latter case we get a new formula of Sylvester type
for the Chow form of any   rational normal scroll.
Here is an explicit example.

\proclaim Example 2.5. (The Chow form of a Hirzebruch surface) \ \
\rm Consider the quadrangle
$$ P \quad = \quad \biggl\{ \,
(m_1,m_2) \in {\bf R}^2 \,: \,
\pmatrix{
\phantom{-}0 & \phantom{-}1 \cr
\phantom{-}1 & \phantom{-}2 \cr
\phantom{-}0 & -1 \cr
-1 & \phantom{-}0 \cr} 
\pmatrix{ m_1 \cr m_2 } \,\leq\,
\pmatrix{ 1 \cr 3 \cr 0 \cr 0 \cr}
\biggr\} \,. $$
The corresponding toric surface is the
rational normal scroll $S_{1,3} $; cf.~[Ha, Example 8.17].
Let $\beta $ be the divisor on $S_{1,3}$ defined by $P$.
Consider three generic elements of $K[x_1,\!\ldots,x_4]_\beta$:
$$
\eqalign{
F_0  \quad &= \quad
   a_1  x_1 x_2^3 
 + a_2  x_1 x_2^2 x_4 
 + a_3  x_1 x_2 x_4^2 
 + a_4  x_1 x_4^3
 + a_5  x_2 x_3 
 + a_6  x_3 x_4 \, ,  \cr
F_1  \quad &= \quad
   b_1  x_1 x_2^3 
 + b_2  x_1 x_2^2 x_4 
 + b_3  x_1 x_2 x_4^2 
 + b_4  x_1 x_4^3
 + b_5  x_2 x_3 
 + b_6  x_3 x_4 \, , \cr
F_2 \quad &= \quad
   c_1  x_1 x_2^3 
 + c_2  x_1 x_2^2 x_4 
 + c_3  x_1 x_2 x_4^2 
 + c_4  x_1 x_4^3
 + c_5  x_2 x_3 
 + c_6  x_3 x_4 \, .\cr}  \eqno (2.11) $$
The quadrangle $3P$ has $10$ interior lattice points,
corresponding to the $10$ monomials of critical degree.
The map $\Phi$ in (2.1) is given by the following
$10 \times 10$-matrix:
$$ \bordermatrix{
   & \, x_2^2  & x_2 x_4  &  x_4^2 &
   &  x_2^2  & x_2 x_4  &  x_4^2 &
   &  x_2^2  & x_2 x_4  &  x_4^2 & 1 \cr
x_1  x_2^5  & \, a_1&   0&  0 &&  b_1&    0&    
0&& c_1&   0 &   0& [125] \cr
x_1  x_2^4  x_4  & \, a_2& a_1&  0 &&  b_2&  b_1&    
0&& c_2&  c_1&   0& [126] + 2 [135] \cr
x_1 x_2^3 x_4^2 &\, a_3& a_2& a_1&&  b_3&  b_2&  
b_1&& c_3&  c_2&  c_1& 
      [235] + 2 [136] + 3 [145] \cr
x_1 x_2^2 x_4^3 &\, a_4& a_3& a_2&&  b_4&  b_3& 
b_2&& c_4&  c_3&  c_2& 
       [236] + 2 [245] + 3 [146] \cr
x_1 x_2 x_4^4 &\, 0 & a_4& a_3&&  0  &  b_4&  
b_3&&  0 &  c_4&  c_3& [345] + 2 [246] \cr
x_1 x_4^5 & \,  0 &  0 & a_4&&  0  &   0 &  
b_4&&  0 &   0 &  c_4& [346] \cr
x_2^3 x_3 & \, a_5&  0 &  0 &&  b_5&   0 &   
0 && c_5&   0 &   0& -[156] \cr
x_2^2 x_3 x_4 &\, a_6& a_5&  0 &&  b_6&  b_5&   
0 && c_6&  c_5&   0& -[256] \cr
x_2 x_3 x_4^2  & \,
  0 & a_6& a_5&&  0  &  b_6&  b_5&&  0 &  
c_6&  c_5& -[356] \cr
x_3 x_4^3 &  0 &\,  0 & a_6&&  0  &    
0&  b_6&&  0 &   0 &  c_6& -[456] \cr}
$$
The border column lists the monomials of critical degree.
The border row gives the multipliers of $F_0,F_1,F_2$ and $J(F)$.
For the coefficients of the Jacobian
$J(F)$ we use the abbreviation
$$ [i\,j \,k] \quad :=  \quad
{\rm det}
\pmatrix{
a_i & a_j & a_k \cr
b_i & b_j & b_k \cr
c_i & c_j & c_k \cr} 
\qquad 
\hbox{for} \,\,\, 1 \leq i < j < k \leq 6 \, .   $$
The determinant of the above $10 \times 10$-matrix
equals the sparse unmixed resultant of (2.11),
i.e., the Chow form of $S_{1,3}$ 
relative to  the given embedding into $P^5$,
by Corollary~2.4. \ \ $\diamond$

\medskip

We close this section with an alternative proof of 
Theorem~1.4, based on
Theorem~2.2.

\vskip .2cm

\noindent {\bf Alternative Proof of Theorem~1.4: \ }
We assume for simplicity that $\ell =1$. The case $\ell > 1$ 
can be dealt with by showing that the matrix of $\Phi$ has a block
decomposition.
We must show that $ \, \RR \cdot \res^X_F(H) \,$ lies in $\, A\,$
for any $H\in C_\rho$. Let ${\cal U}'$ be the intersection of 
$\cal U$ with the Zariski open set where 
all (non identically zero) maximal minors of $\Phi$ do not vanish. 
For $ u \in {\cal U}'$, the ${\bf C}$-linear map
$\Phi_{u}$ is surjective and we can write
$$ H(x) \quad = \quad \sum_{i=0}^n \lambda_i(u;x)\, F_i(u;x) \, + \,
\theta(u)\,J(F^u) \, ,$$
where $\,\theta \, $ 
depends rationally on $u$. By (1.9) we have
$$ \res^X_{F^u}(H) \quad = \quad \gamma \cdot \theta(u) \, ,$$
where $\gamma$ is a rational constant independent of $H$ and $\Efes$.  
This implies that every maximal minor of $\Phi$
which is not identically zero must involve the last column and 
that $\theta(u)$ is unique.
Thus, it follows from Cramer's rule that  
$\res^X_{F}(H)$ may be written as a rational function with
denominator $M$ for all non-identically zero maximal minors
$M$.  Consequently it may also be written as a rational function
with denominator $\RR$.\ \ $\diamond$

\vskip 1 cm

\noi {\bf \S 3. Residues and resultants in the torus}
\smallskip
\noindent
In this section we apply the results of \S 1 to study the
global residue associated with  $n$ Laurent
polynomials in $n$ variables. 
Let $\Delta_1,\ldots,\Delta_n$ be integral  polytopes in
$\R^n$. We form the Minkowski sum
$\Delta := \Delta_1 + \cdots + \Delta_n$ and we 
consider its irredundant presentation
$$\Delta \quad =\quad \{\,m\in\R^n : \ip m {\eta_i} + a_i 
\geq 0\,;\, i=
1,\ldots, s\,\} \,, \eqno{(3.1)}$$
where, as in (1.1), the $\eta_i $
 are the first integral vectors in the inner
 normals to the facets of $\Delta$.  Writing
 $\,a_i^j   =\,  - \, \min_{m\in \Delta_j} 
\langle m, \eta_i \rangle  $,
we get a (generally redundant) inequality presentation
$$\Delta_j \quad = \quad \{\,m\in\R^n : 
\ip m {\eta_i} + a_i^j \geq 0\,;\, i=
1,\ldots, s\,\} \quad
\hbox{for all $j=1,\ldots,n$} \, . $$
The facet normal $\eta_i$ of $\Delta$ supports
a (generally lower-dimensional) face of  $\Delta_j$:
$$ \Delta_j^{\eta_i}\quad := \quad \{\, m\in\Delta_j\,:\, 
 \langle m, {\eta_i} \rangle  = -a_i^{j}\,\}\,. \eqno{(3.2)}$$ 
Consider Laurent polynomials with indetermined
coefficients and Newton polytopes $\Delta_j$,
$$ f_j \quad = \, \sum_{m\in \Delta_j\cap\Z^n} 
u_{jm}\,\cdot\,t^m\,, \eqno{(3.3)}$$  
and introduce the polynomial ring on their coefficients:
$$ A' \quad := \quad \Q[\,u_{jm}\ ;\ 
j=1\ldots,n\ ;\ m\in \Delta_j\cap\Z^n] \,.$$
The leading form of $f_j$ in the direction $\eta_i$ equals
$$f_j^{\eta_i} \quad := \quad
\sum_{m\in \Delta_j^{\eta_i}} u_{jm}\,\cdot\, t^m \,. \eqno (3.4)$$
Since $ \Delta^{\eta_i}\, = \,  \Delta_1^{\eta_i} + 
\ldots + \Delta_n^{\eta_i}$
is a facet of $\Delta$, we may regard $f_1^{\eta_i}, 
\ldots , f_n^{\eta_i}$
as a system of $n$ polynomial functions on an 
$(n-1)$-dimensional torus. We
define $\RR^{\eta_i}$ to be their resultant 
{\sl relative to the ambient lattice $\Z^n$}. 
More precisely, consider the sparse resultant
$\RR_{\Delta_1^{\eta_i},\ldots, \Delta_n^{\eta_i}}$ 
for the support sets $\Delta_1^{\eta_i}\cap \Z^n, \dots,
\Delta_n^{\eta_i}\cap \Z^n$. This is the unique
irreducible polynomial in $A'$ which 
vanishes whenever $f_1^{\eta_i}, \ldots,
f_n^{\eta_i}$ have a common zero in 
$\left( \C^* \right)^n$.
Let $L_j^{\eta_i} \, := \, {\rm aff}_\Z( \Delta_j^{\eta_i} 
\cap \Z^n)$ be the
affine lattice spanned by the integral points in 
$\Delta_j^{\eta_i}$, and
let $L^{\eta_i}\, = \, 
{\rm aff}_\R (\Delta^{\eta_i}) \, \cap \, \Z^n $
be the restriction of $\Z^n$ to the $i$-th facet 
hyperplane of $\Delta$. The index
$\, \ell_i \, : = \, [\, L^{\eta_i} \, : 
\, L_1^{\eta_i} + \ldots +
L_n^{\eta_i} \, ]\,$ is finite.
We define the $i$-th {\it facet resultant} to be
$$\RR^{\eta_i} \quad := \quad
\left( \RR_{\Delta_1^{\eta_i},\ldots, 
\Delta_n^{\eta_i}}\right)^{\ell_i}
\quad  \,\hbox{for} \quad i=1,\ldots,s \, . \eqno (3.5)$$
We now specialize the coefficients $u_{jm}$ in (3.4) 
to complex numbers
such that
$$\RR^{\eta_i}(u) \quad \not= 
\quad 0 \quad \hbox { for } \, i=1,\ldots,s \, .
\eqno (3.6) $$
By Bernstein's Theorem [GKZ, \S 6.2.D, Thm.~2.8],
the hypothesis (3.6) is equivalent to
$${\rm dim}_\C \left( \C[t_1^{\pm 1}, \ldots, t_n^{\pm 1}] / \langle 
f_1, \ldots, f_n \rangle \right) \quad = \quad 
{\rm MV}(\Delta_1,\ldots, \Delta_n) \,,
\eqno (3.6')$$
where ${\rm MV}( \, \cdots \,) $ denotes the {\it mixed volume}.
Let $V$ be the (finite) set of common zeros of 
$f_1,\ldots,f_n$ in the
torus $ \, T = \left( \C^* \right)^n \, $. 
Given any Laurent polynomial
$\,q\in\C[t_1^{\pm 1}, \ldots,t_n^{\pm 1}]$,  the 
{\it global residue} of the differential form 
$$ \phi_q \quad = \quad \frac q{f_1\cdots f_n}\,
\frac {dt_1}{t_1}\wedge
\cdots\wedge \frac {dt_n}{t_n}\,, \eqno{(3.7)}$$ 
is defined as the sum of the local Grothendieck
residues of $\phi_q$, 
at each of the points in $V$:
$$\res_f^T(q) \quad = \quad \sum_{p \in V} 
{\rm Res}_{p,f} (\phi_q) \, .\eqno{(3.8)}$$
We refer to [GH], [AY], and [T] for 
the classical analytic definition of residues and
to [H],  [Ku] or [SS] for the algebraic definition of the
Grothendieck residue.

Note that
$ \, \res_f^T( J^T_f) \, = \, {\rm MV}(\Delta_1,\ldots, \Delta_n)$,
where $J^T_f$ denotes the {\it affine toric Jacobian\/}
$$J^T(f) \quad :=  \quad
\det \bigl( t_k\,\frac {\partial f_j} {\partial t_k}\bigr)_{1\leq 
j,k \leq n}\, . \eqno (3.9)$$
If all the roots of $f_1,\ldots,f_n$ are simple,
i.e.~if $V$ has cardinality ${\rm MV}(\Delta_1,\ldots, \Delta_n)$, 
then 
$$\res_f^T(q) \quad = \quad \sum_{\xi \in V} \frac {q(\xi)}
{J^T(f)(\xi)} \,. \eqno{(3.10)}$$
We conclude from (3.8) or (3.10) that, for fixed $q \in  
\C[t_1^{\pm 1}, \ldots,t_n^{\pm 1}]$, 
the global residue $\res_f^T(q)$
depends rationally on the coefficients $u$. In particular, for any
$m \in \Z^n, \, \res_f^T(t^m)$ is a rational function in $u$ with
$\Q$-coefficients.

Gel'fond and Khovanskii [GK] give a formula for evaluating that
rational function, provided the Newton polytopes $\Delta_1, \ldots,
\Delta_n$ satisfy the following genericity hypothesis:
$$\forall \, i \in \{1,\ldots,s\} 
\quad \exists \, j \in \{1,\ldots,n\}
\, : \, {\rm dim}(\Delta_j^{\eta_i}) = 0 \, . \eqno (3.11) $$
The Gel'fond-Khovanskii formula implies the 
following result, which appears
also in [Z]:

\proclaim Proposition~3.1.   Suppose the
Newton polytopes $\Delta_1,\ldots,\Delta_n$ satisfy {\rm (3.11)}.
Then, for any  $m\in\Z^n$, the residue
$ \res^T_f(t^m)$ is a Laurent polynomial in the coefficients of
$f_1,\ldots,f_n$.

\smallskip

If (3.11) is violated then $ \res^T_f(t^m)$ is generally not a
Laurent polynomial.
In particular, it is never a non-zero Laurent polynomial 
in the unmixed
case $\Delta_1 = \ldots = \Delta_n, \, n \geq 2$. 

Our aim is to characterize the denominator of $ \res^T_f(t^m)$.
 For each $m \in \Z^n$ we define
$$\mu_i^-(m) \quad := \quad - \min\,
\{0, \ip m {\eta_i} + a_i -1\,\}\quad ;
\quad i=1,\ldots,s\,. \eqno (3.12)$$
Geometrically, $\mu_i^-(m)>0$ if $m$ lies beyond 
the facet $\Delta^{\eta_i}$.
We state the main result of this section:
\medskip

\proclaim  Theorem~3.2. Let $f_1,\ldots,f_n$ be generic polynomials
with Newton polytopes $\Delta_1,\ldots,\Delta_n$. For any
 $m\in\Z^n$, the following expression is a 
polynomial in $A$':
$$ \res^T_f(t^m) \cdot
\prod_{i=1}^s \RR^{\eta_i}(f_1^{\eta_i},\ldots,
f_n^{\eta_i})^{\mu^-_i(m)} \, .$$

\medskip
It is easy to derive Proposition~3.1 from Theorem~3.2: If 
$\Delta_j^{\eta_i}
 \, = \, \{m\}$ in (3.11) then  $\RR^{\eta_i} = u_{jm}$ or
$\RR^{\eta_i} = 1$. In fact, (3.11) holds if and only if
$\RR^{\eta_1} \RR^{\eta_2} \ldots \RR^{\eta_s}$ is a monomial. We
present an example where some facet resultants $\RR^{\eta_i}$ are
monomials and others are not.

\medskip

\noi {\bf Example 3.3.} Let $n=2$ and consider the mixed system
$$f_1(t_1,t_2)\quad = \quad a_0t_1+a_1 t_1 t_2 + a_2 t_2^2 \quad ,
\quad f_2(t_1,t_2)\quad = \quad b_0t_2+b_1 t_1 t_2 + b_2 t_1^2 
\, .$$
The Minkowski sum of their Newton triangles is the pentagon
$$\Delta \quad = \quad \Delta_1 + \Delta_2 = 
\biggl\{ \,
(m_1,m_2) \in {\bf R}^2 \,: \,
\pmatrix{
-1 & \phantom{-}0 \cr
-1 & -1 \cr
\phantom{-}0 & -1 \cr
\phantom{-}2 & \phantom{-}1 \cr
\phantom{-}1 & \phantom{-}2} 
\pmatrix{ m_1 \cr m_2 }  + 
\pmatrix { \phantom{-}3 \cr \phantom{-}4 \cr
\phantom{-}3 \cr -3 \cr -3 \cr} \,\geq\,
\pmatrix{ 0 \cr 0 \cr 0 \cr 0 \cr 0 \cr}
\biggr\} \,. $$
The $\Delta$-homogenizations of the input polynomials are
$$F_1 \quad = \quad \frac {x_1 x_2^2 x_3^2}{x_4^2 x_5} 
\, \cdot \, f_1\left(
\frac{x_4^2 x_5}{x_1x_2} , \frac{x_4x_5^2}{x_2x_3} \right) 
\quad = \quad
a_0 x_2 x_3^2 + a_1 x_3 x_4 x_5^2 + a_2 x_1 x_5^3 \, , $$
$$F_2 \quad = \quad \frac {x_1^2 x_2^2 x_3}{x_4 x_5^2} 
\, \cdot \, f_2\left(
\frac{x_4^2 x_5}{x_1x_2} , \frac{x_4x_5^2}{x_2x_3} \right) 
\quad = \quad
b_0 x_1^2 x_2 + b_1 x_1 x_4^2 x_5 + b_2 x_3 x_4^3 \,.$$
Consider the lattice point $m = (3,3)$, which lies 
beyond three facets
of $\Delta$. The global residue of the corresponding monomial 
$t_1^3 t_2^3$ is equal to
$$ \res^T_f(t_1^3 t_2^3) \quad = \quad 
\frac{a_0 a_1 a_2 b_0 b_1 b_2
+ a_0 a_2^2 b_0 b_2^2 - a_1^3 b_0^2 b_2 - a_0^2 a_2 b_1^3}
{a_2 b_2 (a_1 b_1 - a_2 b_2)^3} \,.$$
The denominator can be derived from Theorem~3.2, since
$\, \mu_1^-(m) = \mu_3^-(m) =1, \mu_2^-(m)=3 , 
\mu_4^-(m)= \mu_5^-(m)=0 \,$
and the five facets resultants are 
$$ \RR^{\eta_1} = b_2 \, , \quad 
\quad \RR^{\eta_2} = a_1b_1 - a_2 b_2  \, , \quad 
\quad \RR^{\eta_3} = a_2 \, , \quad
\RR^{\eta_4} = b_0 \, , \quad \hbox{and} \quad
\RR^{\eta_5} = a_0 \, .  \,\,\, \diamond $$    

\bigskip

We shall develop the proof of Theorem~3.2 in several steps.
We first consider the unmixed case $\,P := \Delta_1 = \cdots 
= \Delta_n$.
Let $P$ be presented  as in (1.1) and  ${\cal L}$
the associated line bundle on $X$.
Fix an integer $k_0 > 0$ such that ${\cal L}^{k_0}$ is very ample.
Consider the
mixed sparse resultant $\RR_{k_0} 
:= \RR_{k_0P,P,\ldots, P}$ associated
with the support sets 
$k_0P\cap\Z^n \! ,P\cap\Z^n \!,\ldots, P\cap\Z^n$.
Thus, $\RR_{k_0}$ coincides with the resultant
associated to the line bundles  $\LL^{k_0}\! ,\LL,\ldots,\LL $.  
In the following formula we evaluate $\RR_{k_0}$
at a special monomial section $ t^m $ of $\LL^{k_0}$
and generic sections of  $\LL,\ldots,\LL $.  
Note that the facet resultants ${\cal R}^{\eta_i}$ are
irreducible if ${\cal L}$ is very ample.

\proclaim Lemma~3.4. 
For any $\,m\in k_0P\cap\Z^n \,$ 
we have the following  identity  in $A'$:
 $$ \RR_{k_0}(t^m,f_1,\ldots,f_n) \quad = \quad
\prod_{i=1}^s \,  \RR^{\eta_i}(f_1^{\eta_i},\ldots,
f_n^{\eta_i})^{\ip 
m {\eta_i} + k_0b_i}\, . $$

\noi{\bf Proof: }  Theorem~1.1 in [PSt] gives the following
identity of rational functions:
$$\RR_{k_0}(f_0,f_1,\ldots,f_n) \quad =  \quad
\biggl(\prod_{\xi\in V(f_1,\ldots,f_n)} \!\!\!\!\!\!
f_0(\xi) \,\, \biggr) \cdot
\prod_{i=1}^s  \RR^{\eta_i}(f_1^{\eta_i},
\ldots,f_n^{\eta_i})^{k_0b_i} \, ,
\eqno (3.13) $$
where $f_0,f_1,\ldots,f_n$ are generic polynomials supported in 
$k_0P,P,\ldots,P$.   On the other hand, 
the same result applied to the
 support sets $\{m\},P\cap\Z^n,\ldots,P\cap\Z^n$ gives
$$\prod_{\xi\in V(f_1,\ldots,f_n)} \! \! \! 
\xi^m \quad = \quad
\prod_{i=1}^s \RR^{\eta_i}
(f_1^{\eta_i},\ldots,f_n^{\eta_i})^{\ip m {\eta_i} } \eqno(3.14)$$
since $\RR_{\{m\},P,\ldots,P}(t^m,f_1,\ldots,f_n)= 1$.
Now combine (3.13) and (3.14) for $f_0 = t^m$. \ \ $\diamond$

\bigskip

For $m \in {\bf Z}^n$ and $1 \leq i \leq s$ we abbreviate
$$\mu^+_i(m) \, := \,\max\{0,\ip m {\eta_i} +  nb_i - 1\} 
\quad \hbox{and} \quad
\mu^-_i(m) \, := \, -\min\{0,\ip m {\eta_i} +  nb_i - 1\} \,. $$
This notation distinguishes the facets of $nP$ visible from $m$
from those not visible from $m$.
The following lemma is the unmixed case of Theorem~3.2.

\medskip

\proclaim Lemma~3.5.   Let $f_1,\ldots,f_n$ be generic polynomials
with support in $P$.  Given $m\in\Z^n$,
$$ \res^T_f(t^m) \cdot
\prod_{i=1}^s \RR^{\eta_i}(f_1^{\eta_i},
\ldots,f_n^{\eta_i})^{\mu^-_i(m)} 
\quad \in \quad A'\,.$$

\noi{\bf Proof:} We denote by $F_1,\ldots,F_n$ 
the generic polynomials
in $S_\beta$ obtained from $f_1,\ldots,f_n$ 
by homogenization as in (1.5).
More precisely, if $\,f_i \, = \sum_{m\in P 
\cap {\bf Z}^n} u_{im} t^m
\,$  then
$$ F_i \,=\, F_i(u;x)  \,\,
= \,\, \sum_{m\in P\cap \Z^n} u_{im}\,
(\prod_{i=1}^s x_i^{\ip m {\eta_i} + b_i})\,. \eqno{(3.15)} $$
It is shown in [CD] that the differential form 
$$ \frac {x^{\mu^+(m)}} {x^{\mu^-(m)}\,F_1\cdots F_n} \cdot \Omega$$
is the meromorphic extension to the toric variety $X$
of the form $\,\phi_{t^m} \,$ on the torus $T$ defined in (3.7).
By Theorem~4 in [CD] (or Lemma~3.6 below), there exist
monomials $x^c$ such that $\,\deg(x^{\mu^-(m) + c}) = k_0 \beta \,$ 
for some (arbitrarily large) positive integer $k_0$.  
Whenever the coefficients of $f_1,\ldots,f_n$ lie in 
the Zariski open set where
none of the facet resultants $\RR^{\eta_i}$ vanishes, 
then $F_1,\ldots,F_n$ have no common zeroes at  infinity.
In this case, $\,\{\,x\in X \,: \,F_1(x)=\cdots = F_n(x)=0\, \} 
\, \subset \, T $
and, as shown in [CCD], [CD], the global residue in the torus of 
$\phi_{t^m}$ may be computed as
$$\res^T_f(t^m) \quad = \quad \res^X_F(x^{\mu^+(m) + c}) \,,$$
where $F$ denotes the $(n+1)$-tuple: $\,F_0=x^{\mu^-(m) + c}\,,\,
F_1,\ldots,\,F_n$.

By  Theorem~1.4, the global residue
$\res^T_f(t^m)$ is a rational function with denominator
 $\RR_{k_0} (x^{\mu^-(m) + c}\,,F_1,\ldots,F_n)$. Lemma~3.4 
implies that
$$\RR_{k_0}(x^{\mu^-(m) + c}\,,F_1,\ldots,F_n) \quad =  \quad
\prod_{i=1}^s \, 
\bigl( \RR^{\eta_i}(f_1^{\eta_i},\ldots,f_n^{\eta_i})\bigr)^
{\mu^-_i(m) + c_i} \, . \eqno (3.16) $$
We conclude that the residue
$\,\res_f^T(t^m)$ may be written as a rational function with
denominator the greatest common divisor of all expressions of 
the form
(3.16), where $c = (c_1,\ldots,c_s)$ runs over all
non-negative integer vectors such that
$\deg(x^{\mu^-(m) + c}) = k_0 \beta$ for some integer $k_0 >0$.  
Since unmixed resultants depend on the coefficients of all
polynomials (e.g.~by~[KSZ, Theorem~5.3]), the
facet resultants
$\RR^{\eta_i}(f_1^{\eta_i},\ldots,f_n^{\eta_i})$ 
are powers of distinct 
irreducible polynomials. The proof of Lemma~3.5 
follows from Lemma~3.6 below.
\ \ $\diamond$

\bigskip

\proclaim Lemma~3.6.  
For any non-negative vector $a \in \N^s$ and any 
$i \in \{1,\ldots,s\}\,$  there exists 
a non-negative vector $c\in \N^s$ such that
$\,c_i = 0 \,$ and $\,\deg(x^{a+c}) = k_0\beta \,$ 
for some $k_0\in\N$.

\medskip
\noi{\bf Proof: }  
Let $u^{(1)},\ldots,u^{(\tau)} \in {\bf Z}^n$ be all the vertices
of the lattice polytope $P$ which lie on the facet
$\,P^{\eta_i} = \{\,m\in P : \ip m {\eta_i} + b_i  = 0 \,\}  $.
Their sum $\, u := u^{(1)} + \cdots + u^{(\tau)} \,$ satisfies
$\, \langle u, \eta_i  \rangle + \tau \cdot b_i = 0 \,$ and
$\, \langle u, \eta_j  \rangle + \tau \cdot b_j \geq 1 \,$ for all 
$j \not= i $. Since $\eta_i$ is primitive, we can
find $m \in {\bf Z}^n$ such that $\langle m,\eta_i\rangle = a_i$.
Let $k_0$ be an integer divisible by $\tau$  such that
$$ c_j \quad := \quad
{k_0 \over \tau} \cdot (\langle u, \eta_j  \rangle + \tau \cdot b_j )
  + \langle m , \eta_j \rangle \, - \, a_j  $$
is non-negative for $j=1,2\ldots,s$.
Then  $c = (c_1,\ldots,c_s)$ has the desired properties.
\ \ $\diamond$

\bigskip

We now prove Theorem~3.2 for mixed systems 
of generic Laurent polynomials.

\medskip

\noi{\bf Proof of Theorem~3.2: } 
We shall assume $\,MV(\Delta_1,\ldots, \Delta_n) > 0$.
Otherwise the residue
$ \res^T_f(t^m) $ is zero and Theorem~3.2 trivially holds.

Let $X = X_\Delta$ be the projective toric
variety associated with $\Delta$.  We consider
the homogenization of
the Laurent polynomial $f_j (t_1,\ldots,t_n)$:
$$F_j(x_1,\ldots,x_s) \quad := \quad
\sum_{m\in \Delta_j\cap\Z^n}\,u_{jm}\,\bigl(
\prod_{i=1}^s \,x_i^{\ip m {\eta_i} + a_i^j}\bigr)\,. $$ 
Note that $F_j(x)$ is  generic of degree 
$\alpha_j := [\sum_{i=1}^s a_i^j D_i]$.  Let 
$\alpha := \alpha_1 +\cdots+ \alpha_n = [\sum_{i=1}^s a_i D_i]$.  
For each $j=1,\ldots,n$, let $Q_j$ be a generic polynomial of
degree $\alpha - \alpha_j$ and set $G_j=F_j\,Q_j$. Given a positive
integer $k_0$,  let $F_0$ be a generic polynomial of degree
$k_0\,\alpha$.  Thus
$F_0,G_1,\ldots,G_n$ are homogeneous polynomials of degrees
$k_0\alpha,\alpha,\ldots,\alpha$. For all choices of complex
coefficients in a Zariski open set, they have no common
roots in $X$.  Given a polynomial $H$ of critical degree
$\,\rho(F) := (k_0+1)\,\alpha - \beta_0 \, $ relative
to the $(n+1)$-tuple $F = (F_0,F_1,\ldots, F_n$),
we can compute the toric residue
$\res_F^X(H)$ and, according to the Global Transformation Law 
[CCD, Theorem~0.1]:
$$\res_F^X(H) \quad = \quad \res_G^X(H\cdot Q_1\cdots Q_n)\,\,;\quad
G = (F_0, G_1,\ldots,G_n)\,.$$
Let $\cal R $ be the $(k_0 \Delta, \Delta, \dots, \Delta)$-resultant.
It follows from Theorem~1.4  that the 
specialization $\RR(F_0,G_1,\ldots,G_n)$ is a denominator for the 
rational function $\res_F^X(H)$.
\smallskip
Let $f_0$ denote the dehomogenization of $F_0$, 
let $q_j$ be the dehomogenization of $Q_j$, and set
$\, g_j := f_j\cdot q_j$ for any $j = 1,\ldots, n$.  
Then, $\RR(F_0,G_1,\ldots,G_n)$
agrees with the sparse resultant 
$\RR(f_0,g_1,\ldots,g_n)$ arising from the support sets
$k_0\Delta\cap\Z^n,\Delta\cap\Z^n,\ldots,\Delta\cap\Z^n$.
Given a subset $J\subseteq \{1,\ldots,n\}$, we  denote 
$\tilde f_j := f_j$ if $j\in J$, and $\tilde f_k := q_k$ if 
$k\not\in
J$.  We let $\tilde\Delta_j$ stand for the Newton polytope of
$\tilde f_j$, i.e. $\tilde\Delta_j = \Delta_j$
if $j\in J$, and 
$$\tilde\Delta_k = \Delta_1+\cdots+\Delta_{k-1}+\Delta_{k+1} 
+\cdots+\Delta_n
\  \hbox{ if}\quad k\not\in J\,.$$  
It follows from the Product Formula for sparse mixed resultants 
[PSt, Proposition~7.1] that
$$\RR(f_0,g_1,\ldots,g_n) \quad = \,\prod_{J\subseteq 
\{1,\ldots,n\}}
\RR^J(f_0,\tilde f_1,\ldots,\tilde f_n) \, ,
\eqno{(3.17)}$$
where $\RR^J$ denotes the sparse mixed resultant associated with
the support sets
$$k_0\Delta\cap\Z^n,
\tilde\Delta_1\cap\Z^n,\ldots,\tilde\Delta_n\cap\Z^n\,
\eqno{(3.18)}$$
 relative to the ambient lattice $\Z^n$ as in (3.5).

We now show that the factor $\RR(f_0,\ldots,f_n)$
corresponding, in (3.17), to $J=\{1,\ldots,n\}$ is already a 
denominator of the rational function $\res^X_F(H)$.  Since this
is a function of the coefficients of $f_0,\ldots,f_n$ only, it
suffices to show that every additional factor in (3.17) must involve
the coefficients of some $q_k$, $k=1,\ldots, n$, i.e. 
if $J\not= \{1,\ldots,n\}$, the polynomial
$\RR^J(f_0,\tilde f_1,\ldots,\tilde f_n)$ has positive degree in
the coefficients of some $q_k$, $k\not\in J$. 
But this is a consequence of our assumption
$MV(\Delta_1,\ldots,\Delta_n) > 0$. Indeed,
according to Lemma~1.2 and Corollary~1.1 of [S2], 
it is enough to show
that the collection of supports
$k_0\Delta\cap\Z^n$, $\tilde\Delta_j\cap\Z^n$, $j\in J$
 contains no  proper {\sl essential} subset.  A  subset
which contains $k_0\Delta\cap\Z^n$ cannot be essential since
$\dim(\Delta) = n$ and the cardinality of the subset is at most
$n$. 
On the other hand, no collection of supports 
$ \tilde\Delta_j\cap\Z^n$ can be essential because 
$MV(\Delta_1,\ldots,\Delta_n) > 0$. 

\smallskip
We now complete the proof of Theorem~3.2
similarly to the proof of Lemma~3.5. The algorithm in [CD]
computes $\res^T_f(t^m)$ as the toric residue
$\res^X_F(x^\mu)$ for appropriate monomials $x^\mu$ and 
$F_0(x) = x^\nu$ of degree $(k_0 + 1)\alpha -\beta_0$
and $k_0\alpha$, respectively, where $k_0$ is a positive integer.
For any such choice of $\mu$ and $\nu$, the specialization 
$$\RR(x^\nu, F_1,\ldots,F_n) \quad =  \quad
\prod_{i=1}^s \RR^{\eta_i}(f_1^{\eta_i},
\ldots,f_n^{\eta_i})^{\nu_i}$$
is a denominator of the rational function $\res^T_f(t^m)$.  
Taking the greatest common divisor
over all possible choices and applying Lemma~3.6
yields the theorem.\ \ $\diamond$

\vfill
\eject

\noindent {\bf References} 
\parindent=.45truein
\bigskip 
\item{[AY]}
I. A. A\u{\i}zenberg and A. P. Yuzhakov, 
{\sl Integral representations and
residues in multidimensional complex analysis},  
Translations of
Mathematical Monographs~{\bf 58}. American Mathematical Society, 1983.
\smallskip
\item{[A]}
B.~ Ang\'eniol, {\it R\'esidus et effectivit\'e},  Unpublished
manuscript, 1983.
\smallskip
\item{[B]} V.~Batyrev, 
{\it Variations of the mixed Hodge structure of
affine hypersurfaces in algebraic tori}, Duke Mathematical Journal
{\bf 69} (1993) 349--409.
\smallskip
\item{[BC]} V.~Batyrev and D.~Cox, {\it On the Hodge structure of
projective hypersurfaces in toric varieties}, Duke J.~Math.~{\bf 75}
(1994) 293--338.
\smallskip
\item{[BH]} W.~Bruns and J.~Herzog, {\sl Cohen--Macaulay Rings},
Cambridge Univ. Press, 1993.
\smallskip
\item{[CCD]} E.~Cattani, D.~Cox, and A.~ Dickenstein, 
{\it Residues in
toric varieties}, Compositio Mathematica, to appear, alg-geom 9506024. 
\smallskip
\item{[CD]} 
E.~Cattani and A.~Dickenstein, 
{\it A global view of residues in the torus}, 
Journal of Pure and Applied Algebra, to appear.
\smallskip
\item{[Ch]}
M.~Chardin, {\it The resultant via a Koszul complex}. 
In 
{\sl ``Computational Algebraic Geometry''\/}
(F.~Eysette, A.~Galligo, eds.),
Proceedings MEGA 92, Progress in Math. {\bf 109}, Birkh\"auser, Boston, 1993, 29--39.
\smallskip
\item{[C1]} D.~Cox, {\it The homogeneous coordinate ring of a toric
variety}, Journal of Algebraic Geometry {\bf 4} (1995) 17--50.
\smallskip
\item{[C2]} 
D.~Cox, {\it Toric  residues}, Arkiv f\"or Matematik {\bf 34} (1996) 73--96.
\smallskip
\item{[F]} W.~Fulton, {\sl Introduction to Toric Varieties},
Princeton Univ.~Press, Princeton, 1993.
\smallskip
\item {[GKZ]}
I. Gel'fand, M. Kapranov and A. Zelevinsky, 
{\sl Discriminants, Resultants
and Multidimensional Determinants}, Birkh\"auser,Boston, 1994.
\smallskip
\item {[GK]}
O.~A.~Gel'fond and A.~G.~Khovanskii, 
{\it Newtonian polyhedrons and Grothendieck residues}, 
{Doklady Mathematics} {\bf 54} (1996) 700--702.
\smallskip
\item {[GH]}
P.~Griffiths and J.~Harris, {\sl Principles of Algebraic
Geometry}, John Wiley \& Sons, New York, 1978.
\smallskip
\item {[H]} 
R.~Hartshorne, {\sl Residues and Duality},
Lecture Notes in Math.~{\bf 20}, Springer, 1966.
\smallskip
\item{[Ha]}  J.~Harris, {\sl Algebraic Geometry},
Springer Graduate Texts, New York, 1992.
\smallskip
\item{[J1]} J. P. Jouanolou, 
{\it Singularit\'es rationnelles du r\'esultant},
In Lecture Notes in Mathematics, Vol. 732, Springer, 
1978, pp.~183--213.
\smallskip
\item{[J2]} J. P. Jouanolou, {\it Le Formalisme du R\'esultant},
Advances in Mathematics {\bf 90} (1991) 117--263.
\smallskip
\item{[J3]} J. P. Jouanolou, 
{\it R\'esultant, intersections compl\`etes
et r\'esidu de Grothendieck}, Notes incompl\`etes, 
U. Strasbourg, 1994.
\smallskip
\item{[KSZ]} M.M.~Kapranov, B.~Sturmfels and A.V.~Zelevinsky,
{\it Chow polytopes and general resultants},
Duke Mathematical Journal {\bf 67} (1992) 189--218.
\smallskip
\item{[K1]} A.~G.~Khovanskii, 
{\it Newton polyhedra and toroidal varieties}, 
Functional Analysis and its Applications
{\bf 11} (1977) 289--296.
\smallskip
\item{[K2]} A.~G.~Khovanskii, {\it Newton polyhedra and the 
Euler-Jacobi formula}, 
Russian Math. Surveys~{\bf 33} (1978) 237--238.
\smallskip
\item{[Ku]} E.~Kunz, {\sl K\"ahler Differentials}, Advanced Lectures in 
Mathematics, Vieweg, 1986.
\smallskip
\item{[M]} F.S. Macaulay, 
{\sl The Algebraic Theory of Modular Systems}, 
Cambridge University Press, 1916.
\smallskip
\item{[O]} T.~Oda, {\sl Convex Bodies and Algebraic Geometry},
Springer-Verlag, 1988.
\smallskip
\item{[PSt]}
P. Pedersen and B. Sturmfels, 
{\it Product formulas for resultants and Chow forms}, 
Mathematische Zeitschrift {\bf 214} (1993) 377--396.
\smallskip
\item {[PS]} 
C.~Peters and J.~Steenbrink, {\it Infinitesimal variation of Hodge
structure and the generic Torelli theorem for projective
hypersurfaces\/}. In {\sl Classification of Algebraic and Analytic
Manifolds} (K.~Ueno, ed.), Progress in Math.~{\bf 39},
Birkh\"auser, Boston, 1983, 399--463.
\smallskip
\item{[SS]}  
G.~Scheja and U.~Storch, {\it \"Uber Spurfunktionen bei
vollst\"andigen Durchschnitten}, J.~Reine u.~Angewandte Mathematik 
{\bf 278/9} (1975) 174--190.
\smallskip
\item{[Sta]}  R.~Stanley, {\it Enumerative Combinatorics}, Volume 1,
Wadsworth \& Brooks, Monterey, 1986.
\item{[S1]} B.~Sturmfels, {\it Sparse Elimination Theory}.  In 
{\sl Computational Algebraic Geometry and
Commutative Algebra\/} 
(D.~Eisenbud, L.~Robbiano, eds.), Proceedings, Cortona, June 1991.  
Cambridge University
Press, 1993.
\smallskip
\item{[S2]} B. Sturmfels, 
{\it On the Newton polytope of the resultant}, 
Journal of Algebraic Combinatorics {\bf 3} (1994) 207--236.
\smallskip
\item{[T]} A.~Tsikh, {\sl Multidimensional Residues and Their
Applications}, AMS, Providence, 1992.
\smallskip
\item{[WZ]} J.~Weyman and A.~Zelevinsky,
{\it Determinantal formulas for multigraded resultants},
Journal of Algebraic Geometry {\bf 3} (1994) 569--597.
\smallskip
\item{[Z]} H.~Zhang,  
{\it Sur certains calculs
de r\'esidus pour des syst\`emes de polyn\^omes 
de Laurent interpret\'es
dans le cadre des vari\'et\'es toriques},  
Th\`ese, Universit\'e de Bordeaux, 1996. 

\bye